# Toward a safe supply chain: Incorporating accident, physical, psychosocial, and mental overload risks into supply chain network

***Sajad Karimi, Zaniar Ardalan***

1. Department of System Science and Industrial Engineering, State University of New York at Binghamton, USA

2. Department of Systems and Industrial Engineering, Engineering School, University of Arizona, USA

**Abstract**

Considering health and safety factors in supply chain network design brings safer place for employer and help firm to have better image in the society. There are many health and safety factors overlooked by literature studies of supply chain. This paper takes advantage of the results of occupational safety and health in the transport sector studies and connect this field of science with the supply chain network design. This study incorporates health and safety factors such as accident, physical, psychosocial and mental overload risks as an objective function beside cost and environmental oriented objective functions. We formulated a multi-objective closed loop supply chain network as a mixed integer linear programming model and customized augmented epsilon-constraint algorithm to solve our multi-objective problem to offer multiple choices for decision makers. Eventually, we analyzed the effects of incorporating health and safety factors in supply chain and demonstrated how it will minimize the health and safety risks of supply chain employers, environmental pollution, and the total cost of the network simultaneously.

## 1. Introduction

Network design is one of the key decisions which have a huge effect on the performance of the supply chain network. There are multiple decisions that should be made in this regard such as 1. determining the number and location of facilities like factories, warehouses, distribution centers, and retailers. 2. Allocation of those facilities to each other. 3. Choosing the transportation mode. This decision is made once we establish a supply chain network and it is very costly and time consuming to change the network architecture after establishing it. Therefore, supply chain network should be designed very meticulously.

Nowadays, along with the cost, there are many other objectives that each supply chain network need to satisfy. For example, environmental impacts, social aspects, health and safety of the employers and workers, and so on. In this way, activists, non-governmental organizations and media are demanding firms to take the responsibility of not only their own actions but also other partners in the supply chain [1]. Furthermore, there are legislation enacted by the government to preserve the environment from greenhouse gases in European Union, Australia, and Canada and hopefully, we would have similar legislation all over the world in the near future. This makes firms more concerned about their supply chain design and urge them toward sustainable supply chain management. Furthermore, considering those factors will help firms to have better image in public.

The employer's safety, ergonomic, and health factors in terms of mental, physical, psychological aspects are becoming more and more important, and many researchers in different fields have studied such aspects ([2] and [3]). Based on the study conducted by Schneider and Irastorza [4] published in European Agency for Safety and Health at Work, employer activities can also be categorized based on distance driven, the number of nights away from home, day/night work, type of activity performed. They also stated that drivers are exposed to different health and safety risks [4] such as accident, physical, psychosocial and mental overload risks. These risks are different for each employer based on the activity and nature of the job they perform. The report of European agency for safety and health at work shows that the risk of these health and safety aspects is different for various transportation modes such as railroads, aviation, water and highway transportation sections. Furthermore, Ian Savage [5] analyzed the transportation fatality risk (physical and accident risks) in the United States. He compares the relative risks of the different modes of transportation based on data from 2000 to 2009. Based on the above mentioned research, we can apply the results of different researches in different disciplines and integrate them in the supply chain networks in order to design a desirable low risk network.

Therefore, this paper addresses health and safety risks of supply chain employers as an objective function beside cost and environmental oriented objective functions. We assumed different risks possibilities for different transportation modes. Our network is multi-objective, multi-echelon supply chain in which we considered different transportation modes. This network can be distinguished from the previous studies in the following directions. Firstly, the different health and safety dimensions derived from a comprehensive study in transportation sector are quantified and embedded as a distinct objective besides total costs and environmental impacts. Secondly, the forward and reverse logistics are integrated in a general closed loop supply chain. Thirdly, we modeled our network as a mixed integer linear programming formulation and for tackling the multi-objective nature of the model we modified the augmented Epsilon-constraint method which find the Pareto solutions effectively.

The rest of the paper is organized as follows. The literature reviewed is conducted in Section 2. In Section 3, the problem description and mathematical formulation is presented. In Section 4. we described Augmented Epsilon Constraint. The experimental evaluations are conducted to evaluate the results obtained by the proposed model in Section 5. Finally, Section 6 presents conclusions and future research directions.

## 2. Literature review

In recent years' sustainable supply chain attracted more attention from researchers and practitioners. Tang and Zhou [6] have studied applications of optimization and management sciences in the domain of environmentally and socially sustainable operations. Brandenburg et al. [7] have studied quantitative models that address sustainability aspects in the forward supply chain. There are several papers in the literature studying the Network Design Problem in supply chain. Due to the increasing importance of sustainability, some papers have considered environmental and/or social impacts as additional objective(s) in network design problem. In this case, the problem becomes multi-objective. A fuzzy mathematical programming is developed by Govindan et al [8] for a sustainable multi echelon multi period supply chain network under uncertainty by considering social and environmental objectives as well as cost aspects. The objective function related to the social aspects is considered as increasing the employment opportunity and work harms. In addition, a customized multi-objective particle swarm optimization is developed to find the Pareto solutions. Arampantzi et al [9] presented a new multi objective mixed integer linear programming optimization model for a forward sustainable supply chain with multi period, multi product and multi transportation mode. Three dimensions of sustainability are considered in this research: cost (minimizing investment and operational costs), environmental (reducing greenhouse gases and waste generation) and social (increasing job opportunities, prioritizing societal community development and improving worker conditions). Two solution methods, goal programming and epsilon constraint, are applied to solve the problem. As other works that studied sustainability, sustainability we can mention Dehghanian and Mansour [10] ; Elhedhli and Merrick [11]; Fonseca et al. [12]; Govindan et al. [13]; Kannan et al. [14]; Pishvaee et al. [15]; Wang, Lai and Shi [16], Ardalan et al. [33], and Saeedi et al. [34]. These papers mostly utilized mixed integer programming (MIP) models to formulate network design problems. These models range from simple single objective forward facility location models (Jayaraman and Pirkul [17]) to complex multi-objective closed-loop models (Chaabane et al. [18]). Elhedhli and Merrick [11] considered emission costs alongside fixed and variable location and production costs in network design problem of the forward supply chain. They use a concave function to model the relationship between CO2 emissions and vehicle weight. Wang et al. (2011a) considered the environmental concerns in network design problem of forward supply chain by proposing a multi-objective optimization model that captures the trade-offs between the total costs and the environmental impacts. Pishvaee et al. [15] proposed a credibility-based fuzzy mathematical model for network design problem of Forward supply chain with three stages. Their model aims to minimize both the environmental impacts and the total costs. They show the applicability of the model as well as the usefulness of the solution method in an industrial case study. Fathollahi-Fard et al [19] proposed a multi-objective stochastic closed-loop supply chain and considered social aspects and downside risk as another objectives beside economic objective. For the social aspects, they considered job opportunity and work injuries. In this paper, there is no different transportation modes in the SCM network. For solving the problem, three proposed hybrid meta heuristic algorithms were applied. In another similar research, Hajiaghaei-Keshteli et al [20] worked on a sustainable closed loop supply chain by considering discount supposition for the transportation costs. They developed a mixed integer nonlinear programming model to formulate the three-objective problem in order to minimize the cost, environmental and social aspects.

To avoid the sub-optimality of separate forward and reverse networks, many researchers have integrated forward and reverse networks known as Closed-Loop Supply Chain (Soleimani et al., [21]). Fleischmann et al. [22] suggested a mixed integer linear programming (MILP) model for network design problem of Closed-Loop Supply Chain. Schultmann et al. [23] considered end-of-life vehicle treatment in Germany. They concentrated on the flow of used products and reintegrate reverse flow of used products into their genuine supply chains. Wang and Hsu [24] integrate environmental issues in an integer Closed-Loop Supply Chain model and also develop a genetic algorithm based on a spanning tree structure to solve the problem.

## 3. Problem formulation

In this paper, we introduce a multiple objective mixed integer linear programming (MOMILP) formulation for network design problem of closed-loop supply chain with six echelons to minimize three conflicting objectives: cost of network (fixed cost and transportation costs), negative environmental impacts (greenhouse gas emissions) and negative health risks (accident, psychosocial, physical and mental overload risks). We aim to select the best locations among all potential locations to establish the facilities, best allocation of the facilities to each other and best transportation mode for transporting the products to fulfil our defined objectives. As illustrated in figure 1, the products are produced in manufacturing centers and transported to the distribution centers. Then, the products are delivered to the customers from the distribution centers. We suppose that a fraction of products is not acceptable by customers. Therefore, the customers return the products to the inspection center. The inspection center makes decision about the returned products. If the product is usable, it will be sent to the remanufacturing centers (in the same place of manufacturing centers) to repair and reuse it. If the product is not repairable, the inspection center sends the product to recycle or disposal center. Since our model is multi-mode transportation, all the transportation between facilities can be done by one of the available transportation modes. The proposed model emphasizes the treatment process by considering three types of facilities called treatment centers in the reverse network:

i) Remanufacturing: returned products are remanufactured and/or prepared for reuse;
ii) Recycling: the returned products are recycled and they are used for manufacturing new products; and
iii) Disposal: the returned products are completely disposed as their quality is too low for manufacturing.

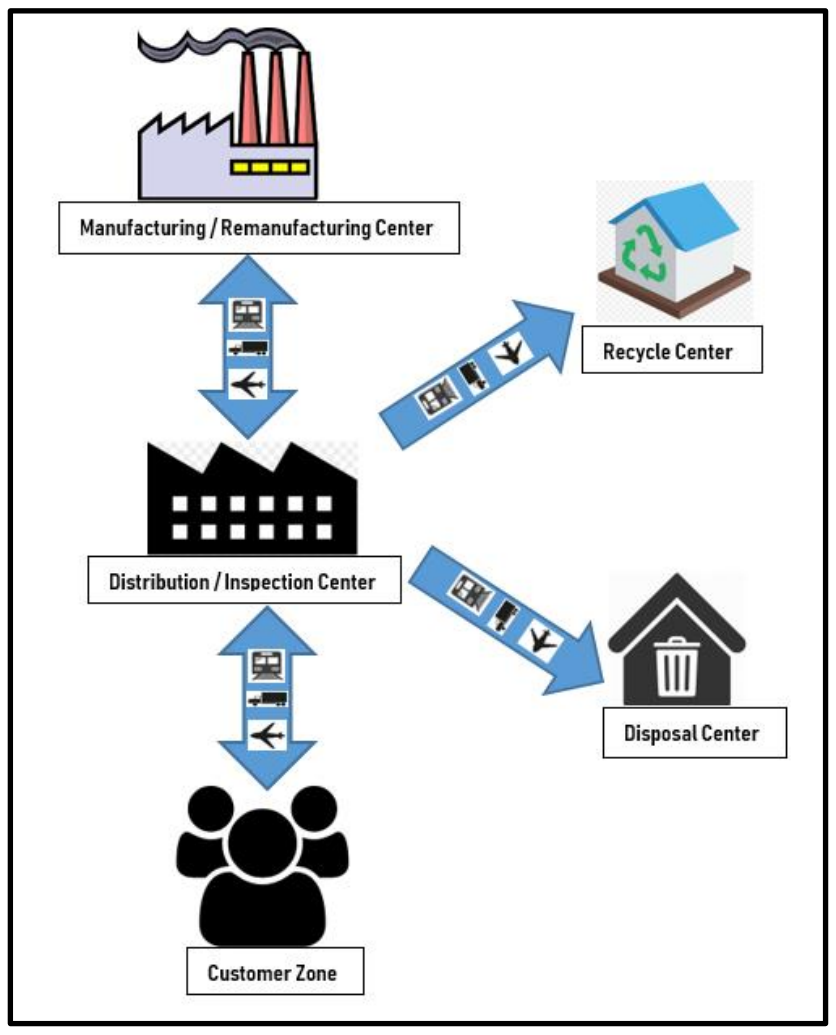

**Figure 1**. The network of the reverse supply chain

Establishing the facilities, processing on products (manufacturing, recycling, etc.) and transportation impose cost to the network and have negative environmental impacts (make greenhouse gas emission). On the other hand, the workers in these activities are encountered with health risks such as accident, physical, psychosocial and mental overload. Different types of activities have different impacts on three aforementioned aspects (economic, environmental and healthcare). For example, the cost of transporting one unit of product with air transportation is higher than road transportation, but the health risks (accident, physical, etc.) in air transportation is much less than road. Therefore, there is a trade-of between three objectives. Our goal is providing some Pareto solutions with reasonable values for all objective functions for the decision makers and they can select one of the solutions based on their policies. To recapitulate, the three objective functions considered in our proposed model are as follows:

**Economic objective**: Minimizing the total fixed cost of establishing facilities and variable costs of transportation.
**Environmental objective**: Minimizing the greenhouse gas emission.
**Health and safety risks objective:** Minimizing the four kinds of risks include accidents, physical, psychosocial and mental overload.
It is noteworthy that the cost, environmental and safety risks for operation activities in the facilities (manufacturing, recycling, etc.) are not considered because we assume that these factors for the operations in all potential locations of each type of facility are equal and do not change the solution, which is a reasonable assumption.

The proposed model is based on the following common assumptions in the literature (Syarif et al [25]; Wang and Hsu [24]; Yao and Hsu [26]). First, all demand of all customers must be satisfied. Second, there is no flow between the facilities of the same echelon.

As a special characteristic of closed-loop logistics and suggested by Laan et al. [27], we assume that the number of the products returned to the inspection centers is a fraction of the customers' demands. In addition, they are allocated to different treatment facilities based on their qualities. The problem under study can be stated as follows.

Given:

- The set of potential sites for locating facilities;
- The set of available transportation modes;
- The demand of customers;
- The cost of opening facilities and shipping materials;
- The environmental impact of opening facilities and shipping materials;
- The fraction of products in each customer zone which is returned to the respective inspection center;
- The fraction of products classified in the inspection centers for each treatment center;
- The capacity of suppliers

Determine:

- The supply chain architecture;
- The selected transportation modes;
- The amount of products to be manufactured at each center;
- The assignment of customers to distribution centers and inspection centers;
- The flow of materials.

The major contributions that distinguish this study from the previous ones are as follows.

- The trade-offs among the three pillars of sustainability is captured by proposing a tri-objectives model.
- The forward and reverse logistics are integrated.
- Different environmental impacts coming from opening facilities, manufacturing at plants and handling products at facilities as well as products damages are modeled.
- The 4 different health and safety factors are also quantified as a separate objective: The number of employers encountered with accident, psychosocial, physical and mental overload risks.

The notations presented in Table 1 are used in the proposed model.

**Table 1.** Parameters of the model.

| Parameter | Description |
|---|---|
| $i,j,k,r,s,t$ | Indexes for manufacturing/remanufacturing, distribution/inspection, customer, recycle, disposal centers and transportation facilities, respectively. |
| $FP_i, FW_j, FR_r, FD_s$ | The fixed cost of opening manufacturing/remanufacturing, distribution/inspection, recycle and disposal centers, respectively. |
| $FPW_{ij}^t, FWC_{jk}^t, FCW_{kj}^t,$ $FWP_{ji}^t, FWR_{jr}^t, FWD_{js}^t$ | The variable cost of transportation between facilities. |
| $EP_i, EW_j, ER_r, ED_s$ | The fixed environmental emission of opening manufacturing/remanufacturing, distribution/inspection, recycle and disposal centers, respectively.. |
| $EPW_{ij}^t, EWC_{jk}^t, ECW_{kj}^t,$ $EWP_{ji}^t, EWR_{jr}^t, EWD_{js}^t$ | The variable environmental emission of transportation between facilities. |
| $\theta_1, \theta_2, \theta_3, \theta_4$ | The weight of health and safety aspect includes accident risks, psychosocial risks, physical risks and risks of mental overload, respectively. |
| $\mu f_i, \pi f_i, \rho f_i, \varphi f_i$ | The number of people encountered with health and safety risks (accident, psychosocial, physical and of mental overload risks, respectively) for establishing a manufacturing / remanufacturing center |
| $\mu w_j, \pi w_j, \rho w_j, \varphi w_j$ | The number of people encountered with health and safety risks (accident, psychosocial, physical and of mental overload risks, respectively) for establishing a distribution / inspection center |
| $\mu c_r, \pi c_r, \rho c_r, \varphi c_r$ | The number of people encountered with health and safety risks (accident, psychosocial, physical and of mental overload risks, respectively) for establishing a recycle center |
| $\mu d_s, \pi d_s, \rho d_s, \varphi d_s$ | The number of people encountered with health and safety risks (accident, psychosocial, physical and of mental overload risks, respectively) for establishing a disposal center |
| $\mu pw_{ij}^t, \pi pw_{ij}^t,$ $\rho pw_{ij}^t, \varphi pw_{ij}^t$ | The number of people encountered with health and safety risks for transporting one unit of product from manufacturing center *i* to distribution center *j* by transportation mode *t* |
| $\mu wc_{jk}^t, \pi wc_{jk}^t,$ $\rho wc_{jk}^t, \varphi wc_{jk}^t$ | The number of people encountered with health and safety risks for transporting one unit of product from distribution center *j* to customer *k* by transportation mode *t* |
| $\mu cw_{kj}^t, \pi cw_{kj}^t,$ $\rho cw_{kj}^t, \varphi cw_{kj}^t$ | The number of people encountered with health and safety for transporting one unit of product from customer *k* to distribution center *j* by transportation mode *t* |
| $\mu wp_{ji}^t, \pi wp_{ji}^t,$ $\rho wp_{ji}^t, \varphi wp_{ji}^t$ | The number of people encountered with health and safety risks for transporting one unit of product from inspection center *j* to remanufacturing center *i* by transportation mode *t* |
| $\mu wr_{jr}^t, \pi wr_{jr}^t,$ $\rho wr_{jr}^t, \varphi wr_{jr}^t$ | The number of people encountered with health and safety risks for transporting one unit of product from inspection center *j* to recycle center *j* by transportation mode *t* |
| $\mu wd_{js}^{t}, \pi wd_{js}^{t}, \rho wd_{js}^{t}, \varphi wd_{js}^{t}$ | The number of people encountered with health and safety risks for transporting one unit of product from inspection center *j* to disposal center *j* by transportation mode *t* |
| $CPF_i, CPR_i, CWF_j, CWR_j, CRC_r, CDS_s$ | The capacity of manufacturing, remanufacturing, distribution, inspection, recycle, disposal centers, respectively. |
| $\alpha$ | The fraction of used products returned from customers |
| $\beta, \delta$ | The fraction of reusable, recyclable of used products returned from customers, inspected from inspection centers. |

*Table 2 presents the decision variables used in the model.*

***Table 2.*** *The decision variables of the model.*

| *Parameter* | *Description* |
|---|---|
| $P_i$ | Binary variables taking value 1 if *i^th^* manufacturing/remanufacturing center is opened; and 0 otherwise. |
| $W_j$ | Binary variables taking value 1 if *j^th^* distribution/inspection center is opened; and 0 otherwise. |
| $C_r$ | Binary variables taking value 1 if *r^th^* recycle center is opened; and 0 otherwise. |
| $D_s$ | Binary variables taking value 1 if *s^th^* disposal center is opened; and 0 otherwise. |
| $X_{ij}^{t}$ | Continues variables for the number of transported products between *i^th^* manufacturing and *j^th^* distribution centers. |
| $Y_{jk}^{t}$ | Continues variables for the number of transported products between *j^th^* distribution center and *k^th^* customer with transportation mode *t* |
| $Z_{kj}^{t}$ | Continues variables for the number of transported products between *k^th^* customer and *j^th^* inspection center with transportation mode *t* |
| $RM_{ji}^{t}$ | Continues variables for the number of transported products between *j^th^* inspection and *i^th^* remanufacturing centers with transportation mode *t* |
| $RC_{jr}^{t}$ | Continues variables for the number of transported products between *j^th^* inspection and *r^th^* recycle centers with transportation mode *t* |
| $DS_{js}^{t}$ | Continues variables for the number of transported products between *j^th^* inspection and *s^th^* disposal centers with transportation mode *t* |

The problem under consideration can be formulated as follows:

$$MIN\ obj_{economical} = \sum_i FP_i.P_i + \sum_j FW_j.W_j + \sum_r FR_r.C_r + \sum_s FD_s.D_s + \sum_i\sum_j\sum_t FPW_{ij}^t.X_{ij}^t + \sum_j\sum_k\sum_t FWC_{jk}^t.Y_{jk}^t + \sum_k\sum_j\sum_t FCW_{kj}^t.Z_{kj}^t + \sum_j\sum_i\sum_t FWP_{ji}^t.RM_{ji}^t + \sum_j\sum_r\sum_t FWR_{jr}^t.RC_{jr}^t + \sum_j\sum_s\sum_t FWD_{js}^t.DS_{js}^t$$

$$MIN\ obj_{environmental} = \sum_i EP_i.P_i + \sum_j EW_j.W_j + \sum_r ER_r.C_r + \sum_s ED_s.D_s + \sum_i\sum_j\sum_t EPW_{ij}^t.X_{ij}^t + \sum_j\sum_k\sum_t EWC_{jk}^t.Y_{jk}^t + \sum_k\sum_j\sum_t ECW_{kj}^t.Z_{kj}^t + \sum_j\sum_i\sum_t EWP_{ji}^t.RM_{ji}^t + \sum_j\sum_r\sum_t EWR_{jr}^t.RC_{jr}^t + \sum_j\sum_s\sum_t EWD_{js}^t.DS_{js}^t$$

$$MIN\ obj_{social} = \sum_i (\theta_1.\mu f_i + \theta_2.\pi f_i + \theta_3.\rho f_i + \theta_4.\varphi f_i).P_i + \sum_j \left(\theta_1.\mu w_j + \theta_2.\pi w_j + \theta_3.\rho w_j + \theta_4.\varphi w_j\right).W_j + \sum_r (\theta_1.\mu c_r + \theta_2.\pi c_r + \theta_3.\rho c_r + \theta_4.\varphi c_r).C_r + \sum_s (\theta_1.\mu d_s + \theta_2.\pi d_s + \theta_3.\rho d_s + \theta_4.\varphi d_s).D_s$$

$$+\sum_{i}\sum_{j}\sum_{t}(\theta_1.\mu pw_{ij}^t+\theta_2.\pi pw_{ij}^t+\theta_3.\rho pw_{ij}^t+\theta_4.\varphi pw_{ij}^t).X_{ij}^t+\sum_{j}\sum_{k}\sum_{t}(\theta_1.\mu wc_{jk}^t$$

$$+\theta_2.\pi wc_{jk}^t+\theta_3.\rho wc_{jk}^t+\theta_4.\varphi wc_{jk}^t).Y_{jk}^t$$

$$+\sum_{k}\sum_{j}\sum_{t}(\theta_1.\mu cw_{kj}^t+\theta_2.\pi cw_{kj}^t+\theta_3.\rho cw_{kj}^t+\theta_4.\varphi cw_{kj}^t).Z_{kj}^t$$

$$+\sum_{j}\sum_{i}\sum_{t}(\theta_1.\mu wp_{ji}^t+\theta_2.\pi wp_{ji}^t+\theta_3.\rho wp_{ji}^t+\theta_4.\varphi wp_{ji}^t).RM_{ji}^t$$

$$+\sum_{j}\sum_{r}\sum_{t}(\theta_1.\mu wr_{jr}^t+\theta_2.\pi wr_{jr}^t+\theta_3.\rho wr_{jr}^t+\theta_4.\varphi wr_{jr}^t).RC_{jr}^t$$

$$+\sum_{j}\sum_{s}\sum_{t}(\theta_1.\mu wd_{js}^t+\theta_2.\pi wd_{js}^t+\theta_3.\rho wd_{js}^t+\theta_4.\varphi wd_{js}^t).DS_{js}^t$$

*Subject to*

$$\sum_{t}\sum_{j}Y_{jk}^t \geq D_k \qquad \forall k \qquad (1)$$

$$\sum_{t}\sum_{i}X_{ij}^t \geq \sum_{t}\sum_{k}Y_{jk}^t \qquad \forall j \qquad (2)$$

$$\sum_{t}\sum_{j}X_{ij}^t \leq CPF_i.P_i \qquad \forall i \qquad (3)$$

$$\sum_{t}\sum_{k}Y_{jk}^t \leq CWF_j.W_j \qquad \forall j \qquad (4)$$

$$\sum_{t}\sum_{j}Z_{kj}^t \geq D_k \times \alpha \qquad \forall k \qquad (5)$$

$$\sum_{t}\sum_{k}Z_{kj}^t \leq CWR_j.W_j \qquad \forall j \qquad (6)$$

$$\sum_{t}\sum_{i}RM_{ji}^t \geq \beta.\sum_{t}\sum_{k}Z_{kj}^t \qquad \forall j \qquad (7)$$

$$\sum_{t}\sum_{j}RM_{ji}^t \leq CPR_i.P_i \qquad \forall i \qquad (8)$$

$$\sum_{t}\sum_{c}RC_{jc}^t \geq \delta.\sum_{t}\sum_{k}Z_{kj}^t \qquad \forall j \qquad (9)$$

$$\sum_{t}\sum_{j}RC_{jr}^t \leq CRC_c.C_r \qquad \forall r \qquad (10)$$

$$\sum_t \sum_d DS_{js}^t \geq (1-\beta-\delta).\sum_t \sum_k Z_{kj}^t \qquad \forall j \qquad (11)$$

$$\sum_t \sum_j DS_{js}^t \leq CDS_s.D_s \qquad \forall s \qquad (12)$$

$$X_{ij}^t, Y_{jk}^t, Z_{kj}^t, RM_{ji}^t, RC_{jr}^t, DS_{js}^t \geq 0 \qquad (13)$$

$$P_i, W_j, C_r, D_s \in \{0,1\} \qquad (14)$$

The first objective is to minimize the total costs of the network. In this objective, the first 4 terms are the fixed costs of establishing the facilities. The 5th to 10th summations are associated with transportation costs. The second objective denotes the environmental impacts (greenhouse gas emissions) of the network. The environmental impacts caused by opening facilities are presented in the first 4 terms. The 5th to 10th terms stand for the environmental impacts of transportation of products. The third objective minimizes the safety aspects of supply chain network. In the first 4 terms the health and safety risks caused by establishing the facilities are presented. The 5th to 10th terms minimizes the health and safety aspects caused by transportation and shipping products between the centers. It is noteworthy that a weight is assigned to each aspect.

Constraint set (1) guarantees that the demand of each customer is satisfied. Constraint set (2) enforces that the total output of each distribution center is less than its total inputs. Constraint sets (3), (4), (6), (8), (10) and (12) ensure that the distributed products from each of the open manufacturing, distribution, inspection, remanufacturing, recycling and disposal centers, respectively, do not exceed their capacity limit. Constraint (5) shows that all the returned products are sent to the inspection centers. Constraint sets (7), (9) and (11) ensure that the returned products of inspection centers are served by the remanufacturing, recycling and disposal centers, respectively. Constraints (13) and (14) define the decision variables.

## 4. Augmented Epsilon Constraint Method

The Epsilon constraint method was presented to solve the multi objective functions by Haimes et. al [28]. The idea behind this method is keeping just one the objectives and consider the other objectives as constraints and restrict them by user defined-values. George Mavrotas [29] proposed an advanced form of this method named Augmented Epsilon Constraint. This method has some advantages compared to the epsilon constraint which are not changing the original feasible region and not depending on the scale of the objectives. In order to use augmented Epsilon constraint, the multi objective problem should be modified as follows:

$$\left(f_1(x) + Epsilon \times (s_2 + s_3 + \cdots + s_n)\right)$$

$$s.t.$$

$$f_2(x) - s_2 = e_2$$

$$f_3(x) - s_3 = e_3$$

...

$$f_n(x) - s_n = e_n$$

$$x \in S \ \ and \quad s_i \in R^+$$

Where “Epsilon” is a small number (between $10^{-3}$ and $10^{-6}$). It is recommended to replace $s_i$ in the above formulation by $\frac{s_i}{r_i}$ to avoid the scaling problems. Where the $r_i$ is the range of $i^{th}$ objective function which can be calculated by payoff table. Therefore, the objective function can be rewritten as follows:

$$\left(f_1(x) + Epsilon \times (s_2/r_2 + s_3/r_3 + \cdots + s_n/r_n)\right)$$

Some different Pareto solutions can be obtained by changing the value of $e_i$ vector. Then, the decision maker can select the most preferred solution among all Pareto solutions. The possible values of $e_i$ vector can be calculated by the following equation.

$$e_i = f_i^{min} - \left(\frac{f_i^{min} - f_i^{max}}{m}\right) \times n, \quad n = 0, 1, 2, \dots, m$$

The above equation divides the range of each constrained objective to m-1 equal intervals and generates m values for the $e_i$. Therefore, the total generated Pareto solutions are $m^p$ in which the p in the number of constrained objectives.

## 5. Numerical example

In this section, in order to evaluate our proposed model, a numerical experiment is designed. In this example, we consider 20 customer zones and the number of potential locations for manufacturing/remanufacturing, distribution/inspection, recycle and disposal centers are 5, 8, 3 and 3 respectively. Also, 3 transportation modes (rail, road and air) are considered to transport the products between the facilities. The demands of the customers are varied and considered between 100 and 150 products per year. The capacity, fixed cost and CO2 emission of opening each facility are selected of the intervals shown in table 3. Note that, although the manufacturing and remanufacturing centers are at the same location, each center has a different capacity for the manufacturing and remanufacturing of the products. In this example, we consider CO2 emissions (as the most primary greenhouse gas) to reflect the negative environmental impacts.

**Table 3**: The capacity, fixed cost and CO2 emission of opening each facility

| Facility | Capacity (Product per year) | Fixed Cost for Opening ($) | CO2 emission (g) |
|---|---|---|---|
| Manufacturing Center | [800,1200] | [400000, 600000] | [4000,8000] |
| Remanufacturing Center | [80,120] | | |
| Distribution Center | [400,600] | [250000, 350000] | [3000,5000] |
| Inspection Center | [80,120] | | |
| Recycle Center | [160, 240] | [150000, 250000] | [2000,4000] |
| Disposal Center | [160, 240] | [150000, 250000] | [2000,4000] |

In this paper we consider the same variable cost, environmental and safety risks of operations for each center. Therefore, this type of parameters do not have any impact on the optimal solution. Thus, we can remove variable operation parameters of the model.

As mentioned before, the products are transported between facilities by 3 transportation modes: rail, road and air. We know that the capacity, travel time and cost of each transportation mode are different. We use the real data of transportation cost in our numerical experiment. Based on the data presented by "Business logistics/supply chain management" [30], the transportation cost for three aforementioned modes are as table 4. In this table, the cost is based on ton-mile for 1995. We inflated the values for 2019 and since the transportation cost in our mathematical modelling is based on each product, we calculated the transportation cost of a product between centers (we supposes that each product is 50 kg and the distances are given, so the transportation cost per product between each pair of centers is easily calculated). Similar to environmental impacts of opening the facilities, we consider CO2 emission for negative environmental impact of transportation. For each mode of transportation, there is a real data of CO2 emissions [31]. Table 4 shows the amount of CO2 (gram) emitted by trucks, trains and airplanes per tonne-mile. As can be seen, the airplane has the most and train has the least CO2 emission. By similar calculation for cost, we calculate the CO2 emission for transporting one unit of product between centers for each transportation mode based on the distance between centers.

**Table 4:** The transportation cost for three aforementioned modes

| Transportation mode | Variable cost ($/ ton-mile) | CO2 emission (g/tonne-km) |
|---|---|---|
| **Rail** | 0.03 | 22 |
| **Road** | 0.25 | 62 |
| **Air** | 0.59 | 602 |

The number of workers associated with the health and safety risks for transporting one unit of product in one mile are as table 5. We extracted the real data for the accident rate of transportation modes (fatality and injury rate) of a report generated by US Department of Transportation [32] and modified the numbers based on one unit of product for our numerical experiment. Based on the distance between nodes, the total worker encountered with health and safety risks are calculated and imported to the mathematical modelling as input parameters. As can be seen, the rate of accidents in the road transportation is much greater than other transportation mode while other health and safety risks (physical, psychosocial and mental overload) are higher in rail transportation. Table 6 shows the health and safety risks of opening the facilities. In this numerical experiment, the weight of accident risk ($\theta_1$)is considered 40% and the weight of other risks ($\theta_2, \theta_3, \theta_4$)are considered 20%. The values of $\alpha, \beta$ and $\delta$ are 0.2, 0.4 and 0.3 respectively.

**Table 5**: Number of workers associated with health and safety risk per unit of product-mile

| Health and Safety Risks | Road (#worker/ product-mile) | Rail (#worker/ product-mile) | Air (#worker/ product-mile) |
|---|---|---|---|
| **Accident** | 0.00005 | 0.000005 | 0.00001 |
| **Psychosocial** | 0.0025 | 0.005 | 0.0005 |
| **Physical** | 0.0025 | 0.005 | 0.0005 |
| **Mental overload** | 0.0025 | 0.005 | 0.0005 |

**Table 6**: Number of workers associated with health and safety risk for opening the facilities

| Health and Safety Risks | Manufacturing / Remanufacturing Center | Distribution / Inspection Center | Recycle Center | Disposal Center |
|---|---|---|---|---|
| **Accident** | 3 | 2 | 2 | 2 |
| **Psychosocial** | 6 | 4 | 4 | 4 |
| **Physical** | 6 | 4 | 4 | 4 |
| **Mental overload** | 4 | 3 | 3 | 3 |

The problem is formulated in the GAMS 22.2 and is run on a PC with 2.1 GHz Intel core i3 processor and 4 GB RAM memory. For the first analysis, we solve the problem based on each objective function separately. For example, in the first trial, the problem is solved based on the first objective function (economic) and the value of other two objective functions are calculated based on the optimal value of decision variables. This process is repeated for the environmental and health and safety objective functions. The value objective functions for the three trials (payoff table) are shown in table 7.

**Table 7**: payoff table of numerical experiment

| | Value of objective functions | | |
|---|---|---|---|
| | Economic | Environmental | Social |
| **Trial 1 (Objective function=Economic)** | **3,319,938** | 50,649 | 454 |
| **Trial 2 (Objective function=Environmental)** | 3,595,870 | **46,145** | 459 |
| **Trial 3 (Objective function=Social)** | 7,045,747 | 104,930 | **124** |

In the first trial with the economic objective function, the rail transportation is selected for all transportation of the network. It is reasonable because rail transportation has the least cost among all transportation modes.

On the other hand, in trial 3 withhealth and safety objective function, airplane is selected for the entire network because the air transportation has the best transportation mode in terms of the health and safety risks. As can be seen, the cost of supply chain in this trial is much higher than the first trial. By similar analysis, CO2 emission is the lowest level in second trial with environmental objective function. Therefore, there is a tradeoff between three conflicting objectives. In this situation, in order to obtain a solution which provides a reasonable value of all three objective functions we use Augmented Epsilon Constraint as a solving method. As mentioned in section 4, augmented epsilon constraint provides some Pareto solutions with different levels of value for each objective function. The manager can select each solution based on the policy of company. In this method we selected 10 cuts and obtained 100 Pareto solutions. One of the solutions is as table 8. As can be seen, even though none of the objectives are not in optimal value, all of them are in a reasonable value. We can see that with considering the health and safety aspects in designing the supply chain network, physical and psychological risks of the workers reduce considerably. In this example, based on one of the Pareto solutions, the number of workers encountered with health issues are decreased by 37%.

**Table 8:** The value of objective functions for one of the Pareto solutions

| **Objective function** | **Value** |
|---|---|
| Economic | 4454895 |
| Environmental | 50198 |
| Social risks | 284 |

## 5. Conclusion

In this paper, we presented a multi-objective optimization model for closed loop six echelon supply chain network design model with multi-mode transportation. The main motivation of the authors is considering safety and health risks besides economic and environmental aspects in designing the supply chain network. Neglecting the health and safety aspects causes many physical and psychological risks which most of them have permanent negative impacts on the workers. Therefore, it is an indispensable task to deal with safety risks in supply chain problems. Thus, three conflicting objectives are considered in our proposed model: economic, environmental, and safety and health. For the economic objective, the model minimizes the total transportation and fixed establish costs and reducing the greenhouse gas emission is considered in environmental objective. In the third objective, minimizing four safety and health risks (accident, psychosocial, physical and mental overload) is considered. A multi-objective mixed integer linear programming is developed to formulate the problem. In order to solve the multi-objective problem and find the pareto solutions, we utilized the augmented epsilon constraint method. To evaluate the results obtained by the proposed model, we conducted a numerical experiment in which a mid-size problem based on the real data was generated and the problem was solved by proposed model in two scenarios. In the first scenario, we solved the problem based on each single objective separately and in second scenario, all three objectives were considered and the problem solved by epsilon constraint method to obtain the pareto solutions.  Analyzing and comparing the solutions revealed that using our proposed multi-objective model provides much reasonable solutions compared with traditional closed-loop supply chain model with single objective (usually cost). In the first scenario (single objective), the value of corresponding objective is optimal, but the others are far from the optimal value and not acceptable in real cases. In contrast, the value of all objectives obtained by our proposed multi-objective model are near optimal and reasonable in real world. Our proposed model provides the solutions that reduce the safety risks significantly by keeping the other objectives in a reasonable value.